\documentclass[11pt,twoside]{amsart}

\usepackage[english]{babel}
\usepackage[utf8]{inputenc}
\usepackage{csquotes}
\usepackage{amsfonts}
\usepackage{hyperref}
\usepackage{etoolbox}
\usepackage{amsmath}
\usepackage{amssymb}
\usepackage{amsthm}
\usepackage{pifont}
\usepackage{caption}
\usepackage{subcaption}
\usepackage{mathrsfs}
\usepackage{graphicx}
\usepackage{enumerate}
\usepackage[all]{xy}
\usepackage{enumitem}
\usepackage{geometry}
\usepackage{amsfonts}
\usepackage{fancyhdr}
\usepackage{calc}
\usepackage{cancel}
\usepackage{abraces}
\usepackage{tikz,tikz-3dplot}
\usetikzlibrary{decorations.markings}
\usetikzlibrary{shapes,arrows,plotmarks}

\usepackage{color}
\definecolor{linkColor}{rgb}{0.0,0.0,0.554}
\definecolor{citeColor}{rgb}{0.0,0.0,0.554}
\definecolor{fileColor}{rgb}{0.0,0.0,0.554}
\definecolor{urlColor}{rgb}{0.0,0.0,0.554}
\definecolor{promptColor}{rgb}{0.0,0.0,0.589}
\definecolor{brkpromptColor}{rgb}{0.589,0.0,0.0}
\definecolor{gapinputColor}{rgb}{0.589,0.0,0.0}
\definecolor{gapoutputColor}{rgb}{0.0,0.0,0.0}
\usepackage{fancyvrb}
\usepackage{adjustbox}
\usepackage{helvet}

\definecolor{cof}{RGB}{219,144,71}
\definecolor{pur}{RGB}{186,146,162}
\definecolor{greeo}{RGB}{91,173,69}
\definecolor{greet}{RGB}{52,111,72}

\tdplotsetmaincoords{70}{165}

\mathchardef\mhyphen="2D 

\title[Extensions of the rational Cherednik algebra and generalized KZ functors.]{Extensions of the rational Cherednik algebra and generalized KZ functors.}

\author{Fallet Henry}
\address{\newline
Université de Picardie,
\newline Département de Mathématiques et LAMFA (UMR 7352 du CNRS),
\newline 33 rue St Leu,
\newline F-80039 Amiens Cedex 1,
\newline France }
\email{henry.fallet@u-picardie.fr}

\theoremstyle{plain}

\newtheorem{thm}{Theoreme}[section]
\newtheorem{prp}[thm]{Proposition}
\newtheorem{dfn}[thm]{Definition}

\newtheorem{lem}[thm]{Lemma}

\newcommand{\fonction}[5]{\begin{array}{l rcl}
#1: & #2 & \longrightarrow & #3 \\
    & #4 & \longmapsto & #5 \end{array}}
\newcommand{\fonctionb}[9]{\begin{array}{l rcl}
#1: & #2 & \longrightarrow & #3 \\
    & #4 & \longmapsto & #5  \\
    & #6 & \longmapsto & #7 \\
    & #8 & \longmapsto & #9 \\
    \end{array}}

\newcommand{\N}{N_{W}(W_{0})}

\newcommand{\C}{\mathbb{C}}
\newcommand{\fa}{\forall \,}
\newcommand{\lra}{\longrightarrow}

\newcommand{\s}{\sigma}

\newcommand{\mcf}{\mathcal{F}}
\newcommand{\mcc}{\mathcal{C}}

\newcommand{\mco}{\mathcal{O}}
\newcommand{\mcd}{\mathcal{D}}
\newcommand{\mca}{\mathcal{A}}
\newcommand{\mcl}{\mathcal{L}}
\newcommand{\mcr}{\mathcal{R}}
\newcommand{\mcg}{\mathcal{G}}

\newcommand{\mcb}{\mathcal{B}}

 \normalbaselines
\makeatletter
\newlength\@SizeOfCirc%
\newcommand{\CircleArrowRight}[1]{%
    \setlength{\@SizeOfCirc}{\maxof{\widthof{#1}}{\heightof{#1}}}%
    \tikz [x=1.0ex,y=1.0ex,line width=.12ex]%
        \draw [->,anchor=center]%
            node (0,0) {#1}%
            (0,0.8\@SizeOfCirc) arc (85:-240:0.8\@SizeOfCirc);%
}%
\newcommand{\CircleArrowLeft}[1]{%
    \setlength{\@SizeOfCirc}{\maxof{\widthof{#1}}{\heightof{#1}}}%
    \tikz [x=1.0ex,y=1.0ex,line width=.12ex]%
        \draw [<-,anchor=center]%
            node (0,0) {#1}%
            (0,0.8\@SizeOfCirc) arc (85:-240:0.8\@SizeOfCirc);%
}%
\makeatother


\date{\today}

\begin{document}

\begin{abstract}
Ginzburg, Guay, Opdam and Rouquier established an equivalence of categories between a quotient category of the category $\mco$ for the rational Cherednik algebra and the category of finite dimension modules of the Hecke algebra of a complex reflection group $W$. 
We establish two generalizations of this result.
On the one hand to the extension of the Hecke algebra associated to the normaliser of a reflection subgroup and on the other hand to the extension of the Hecke algebra by a lattice.
\end{abstract}

\maketitle

\tableofcontents

\setcounter{section}{-1}
\section{Introduction}
In $1967$, Yokonuma \cite{yokonuma1967structure} introduce a generalization of the Hecke algebra in the context of finite reductive group. Let $G$ be a finite Chevalley group. Let $B$ a Borel subgroup of $G$. Let $U$ be the radical unipotent of $B$ and $T$ the maximal torus of $G$. Let $(W,S)$ be the Weyl group of $G$. Let $\textbf{k}$ be a commutative ring. The Yokonuma-Hecke algebra is the algebra of functions from $G$ to $\textbf{k}$ bi-invariant by the action of $U$ endow with a convolution product. 
A subalgebra of the Yokonuma-Hecke algebra of type $A$ was introduced by Aicardi and Juyumaya in 2014 \cite{aicardi2014markov} for the definition of invariants of tied knots, and was therefore named "algebra of braids and ties".


In $2016$, I.Marin \cite{marin2018artin} observed that $BT_{n}$ is an extension of the Hecke algebra  of type $A$ by the lattice of reflexion subgroups of the symmetric group. Furthermore, he defined for any Coxeter system $(W,S)$ an extension $\mcc_{W}$ of the Hecke algebra  $H(W)$ by the lattice of full reflection subgroups. By full reflection subgroup we mean that $W_{0}$ is a reflection subgroup of $W$ and for any reflection  in $W_{0}$, all the reflection with the same reflecting hyperplane belong to $W_{0}$.

Moreover, I.Marin generalized this construction to any finite complex reflection group. Let us first recall the construction of the Hecke algebra associated to a finite complex reflection group.

\medskip

In 1998, Broue, Malle and Rouquier  introduce in their seminal paper \cite{rouquier1998complex}, the Hecke algebra  of a finite complex reflection group. Let $V$ be a   complex $\C$-vector space of finite dimension. Let  $W$ be a finite complex reflection group. Let $\mcr$ be the set of all reflections of $W$ and let $\mca$ be the arrangement of  reflecting hyperplanes of $W$. We denote  by $X$ the complement of $\mca$ in $V$. Let $(k_{H,j})_{j \in \{0, \ldots, m_{H}-1\}}$ be  a set of complex numbers indexed by a reflecting hyperplane of $W$ and integers $j$, where $m_{H}$ is the order of the pointwise stabilizer of $H$ in $W$, such that for all $w \in W$, $k_{w(H),j} = k_{H,j}$ for all integers $j$. Let $B(W):= \pi_{1}(X/W)$ be the braids group of $W$, which is generated by braided reflections $\s_{H}$.
The Hecke algebra of $W$ is $H(W):= \frac{\C B(W)}{I}$ where  $I$ is a two-sided ideal of $\C B(W)$ generated by the relation $\s_{H}^{m_{H}} = \sum_{k=0}^{m_{H}-1}k_{H,j}\s_{H}^{j}$.

Let us now recall the construction made by I.Marin for $W$ a complex reflection group.

\textbf{The algebra} $\mcc(\mcl,W)$.
In, \cite{marin2018artin} and \cite{marin2018lattice}, Marin has defined an extension of $H(W)$ by the lattice $\mcl$ of full reflection subgroups of a finite complex reflection group $W$. We denote by $\C \mcl$ the Möbius algebra of $\mcl$ (\cite{stanley2011enumerative} definition 3.9.1). It is a free $\C$-module generated by $e_{\lambda}$, $\lambda \in \mcl$ with an inner multiplication given by $e_{\lambda}.e_{\mu} = e_{\lambda \vee \mu}$ where $\vee$ is the join rule of $\mcl$. The group $W$ acts by automorphisms on $\C \mcl$ so does $B(W)$. We can consider the algebra $B(W) \ltimes \C \mcl$. The algebra $\mcc(\mcl, W)$ is the quotient of the algebra $\C \mcl \rtimes B(W)$ by the two-sided ideal $\mathcal{J}$ of $\C \mcl \rtimes B(W)$ generated by the relation $\s_{H}^{m_{H}} - 1 = (\sum_{j = 0}^{m_{H}-1}k_{H,j}\s_{H}^{j}+1)e_{\lambda}$ with $\lambda \in \mcl$.  We have a surjective morphism of algebras: $\mcc(\mcl,W) \twoheadrightarrow H(W)$.

Let us remark that for $W_{0} \in \mcl$ then $\text{stab}_{W}(W_{0}) = N_{W}(W_{0})$ the normalizer of $W_{0}$ in $W$. We have a bijection between the orbit $[W_{0}]$ of $W_{0}$ under $W$ and the quotient $W/\N$.

Since $\C \mcl \rtimes W \simeq \C^{\mcl} \rtimes W \simeq \bigoplus_{[W_{0}] \in \mcl/W}\text{Mat}_{[W_{0}]}(\C \N)$  (\cite{marin2018lattice} proposition 2.1), we can expect that $\mcc(\mcl, W)$ is Morita equivalent to an algebra of matrices involving the algebra $\C \N$.  Indeed, the algebra  $\mcc(\mcl,W)$ is Morita equivalent to a direct sum of matrices algebras with coefficient in a new Hecke algebra. The Hecke algebra of $\N$. 

\medskip

\textbf{The Hecke algebra of the pair} $(W_{0},W)$.
Let $W_{0}$ be a reflection subgroup of $W$. We denote $\N$ the normalizer of $W_{0}$ in $W$. Let $\mcr_{0}$ be the set of reflections of $W_{0}$, $\mcr_{0}:= \mcr \cap W_{0}$. Let $\mca_{0}$ be the arrangement of reflecting hyperplanes of $W_{0}$, $\mca_{0}:= \{ Ker(s-1) | s \in \mcr_{0}\}$. Let $X_{0}$ be the complement of the hyperplane arrangement $\mca_{0}$ in $V$.

\medskip

 Let $\hat{B_{0}}$ be the group $\pi_{1}(X/\N)$. Let $J$ be the  two sided-ideal of $\C \widetilde{B_{0}}$ defined by $\s_{H}^{m_{H}} = 1$ where $\s_{H}$ is a braided reflection associated to $H \in \mca \setminus \mca_{0}$ and $\s_{H}^{m_{H}} = \sum_{j=0}^{m_{H}-1}k_{H,j}\s_{H}^{j}$ with $H \in \mca_{0}$. The Hecke algebra of $\N$ is the quotient of $\C \widetilde{B_{0}}$ by $J$. It is denoted $H(W_{0},W)$. There is a surjective morphism of algebras from  $H(W_{0},W)$ to $H(W_{0})$.

Then we have a structure theorem

\begin{thm}(\cite{marin2018lattice} theorem 2.3) \label{structuretheorem}It exists an isomorphism of $\C$-algebras,
$$\mathcal{C}(\mcl, W) \simeq \bigoplus_{[W_{0}] \in \mcl /W}Mat_{|[W_{0}]|}(H(W_{0},W))$$
\end{thm}
\medskip
\textbf{The Knizhnik–Zamolodchikov functor.}
In $2003$, Ginzburg, Guay, Opdam and Rouquier \cite{ginzburg2003category} defined a category $\mco$  associated to a Cherednik algebra $A(W)$ of a finite complex reflection group $W$ \cite{etingof2002symplectic}. This category of modules over $A(W)$ is endowed with finiteness conditions, based on the Bernstein-Gelfand-Gelfand  category $\mco$ for a semisimple Lie algebra. They constructed a functor, called Knizhnik–Zamolodchikov functor  $KZ$, from this category $\mco$ to the category of finite dimension modules over the Hecke algebra of $W$. The construction of this functor made a critical use of the Dunkl-Opdam differential operators.  This functor induces an equivalence of categories  between the category $\mco/\mco_{tor}$ of object in $\mco$ supported outside of $X$ and the category of all finite dimensional representations of the Hecke algebra of a complex reflection group.  
\medskip

\textbf{The aim} of this article is to establish two results of the same kind as in \cite{ginzburg2003category} for the case of the algebras $H(W_{0},W)$ and $\mcc(\mcl,W)$.
We start by constructing a Cherednik algebra for $\C \N$ as a symplectic reflection algebra for $\N$ \cite{etingof2002symplectic} depending on a parameter $t$. In this article, only the $t \neq 0$ case will be investigated and up to renormalization, we can consider $t = 1$. Then we introduce a family of commutative differential operators, $\N$-equivariant by considering the covariant derivative of a flat connection involving a $1$-differential form defined in \cite{gobet2020hecke} (proposition 2.6). In order to prove the commutativity we use arguments from \cite{etingof2010lecture} theorem 2.15. This allows us to define a Dunkl embedding of $A(W_{0}, W)$ inside an algebra of differential operators on $X$ , $\N$-equivariant $\mcd(X)\rtimes \N$. 

\medskip

Then we define an Euler element $eu_{0}$, standard objects,  a highest-weight category $\mco(W_{0}, W)$ and a functor $KZ_{0}$ from $\mco(W_{0},W)$ to the category of finite dimensional representation of $H(W_{0},W)$. From an object $\mco(W_{0},W)$, we get a $\mcd(X)^{\N}$-modules by applying a localization functor, the Dunkl embedding and some classical Morita's equivalence. We use a general result of algebraic geometry \cite{cannings1994differential} theorem 3.7.1 to prove $\mcd(X)^{\N} \simeq \mcd(X/\N)$.

\medskip

We define a flat connection on standard objects of $\mco(W_{0},W)$ and prove that this connection has regular singularities on $V$. Afterwards, we prove that we can apply the Riemann-Hilbert-Deligne equivalence for any object of $\mco(W_{0},W)$. We get a $\C \pi_{1}(X/\N)$-module finitely generated. Furthermore, the \cite{gobet2020hecke} proposition 2.6 implies that the monodromy action of $\C \pi_{1}(X/W)$ factors through an action of $H(W_{0},W)$.  As a result, we had built a functor $$KZ_{0}: \mco((W_{0},W) \rightarrow H(W_{0},W)\text{-mod}_{f.d}$$

\medskip

Since this functor is exact we can prove that this functor is representable by a projective object of $\mco(W_{0},W)$, this allows us to prove that the category image of $KZ_{0}$ is actually the category $H(W_{0},W)\text{-mod}_{f.d}$. By killing all the objects of $\mco(W_{0},W)$ with support outside of $X$ we get
\begin{thm}
The functor $$\underline{KZ_{0}}: \mco(W_{0},W)/\mco_{tor}(W_{0},W) \rightarrow H(W_{0},W)\text{-mod}_{f.d}$$ is an equivalence of categories.
\end{thm}

\medskip

 In the second part of this article, we develop a similar construction in the context of the lattice extension. We start by defining a Cherednik-algebra for the pair $(\mcl,W)$ as an algebra with a triangular decomposition \cite{bellamy2018highest} (definition 3.1) depending on a parameter $t$. We consider the case $t \neq 0$ and up to a renormalization $t = 1$.  Then we introduce a family of commutative differential operators, $W$-equivariant by considering the covariant derivative of a flat connection involving a $1$-differential form defined in \cite{marin2020truncations} proposition 2.5.  In order to prove the commutativity we use arguments from \cite{etingof2010lecture} theorem 2.15. Whereas in the previous case the action of $(\mcd(X) \otimes \C \mcl)\rtimes W$ on $\C (V)$ is not obviously faithful. In order to prove the faithfulness, we prove a more general result:
\begin{prp}
    Let $R$ be a simple ring with unity. Let $G$ be a group of outer automorphisms of $R$. Let $X$ be a finite set. If $G$ acts transitively on $X$, then $(R \otimes \C^{X})\rtimes G$ is simple.
    \end{prp}
  This allows us to define a Dunkl embedding of $A(W_{0}, W)$ inside the algebra $(\mcd(X) \otimes \C \mcl) \rtimes W$. Then we define an Euler element $\widetilde{eu}$, standard objects with $Irr(\C \mcl \rtimes W)$, and a highest weight category $\mco(\mcl,W)$. We prove that we can endow an object of $\mco(\mcl,W)$ with flat connection over a trivial vector bundle on $V$ with fiber $E \in Irr(\C \mcl \rtimes W)$ and we prove that this connection has regular singularities on $V$. This allows us to apply the Riemann-Hilbert-Deligne correspondence, we then  get a finitely generated  $(\C B(W) \otimes \C \mcl)\rtimes W$-modules. Furthermore, \cite{marin2018lattice}, the monodromy action factorizes through an action of $\mcc(\mcl,W)$. Consequently, we have defined a functor: $\widetilde{KZ}: \mco(\mcl,W) \rightarrow  \mcc(\mcl,W)\text{-mod}_{f.d}$.

\medskip  
  
  Since this functor is exact we can prove that this functor is representable by a projective object of $\mco(\mcl,W)$, this allows us to prove that the category image of $\widetilde{KZ}$ is the category $\mcc(\mcl,W)\text{-mod}_{f.d}$.
  Then we  kill all the object of $\mco(\mcl, W)$ with support outside of $X$, and  we get
\begin{thm}
The functor $\underline{\widetilde{KZ}}: \mco(\mcl,W)/\mco_{tor}(\mcl,W) \rightarrow \mcc(\mcl,W)\text{-mod}_{f.d}$ is an equivalence of categories.
\end{thm}

\section{The rational Cherednik algebra of the normalizer \texorpdfstring{$\N$}- and the functor \texorpdfstring{$KZ_{0}$}-.}
\subsection{The rational Cherednik algebra of the normalizer}
Let $\delta = \prod_{H \in \mca}\alpha_{H}$ where $\alpha_{H}\in V^{*}$ such that $Ker(\alpha_{H}) = H$. This is  an element of $\C[V]$ and $\delta$ vanishes on $\bigcup \mca$. 

We define the rational Cherednik algebra  of $\C \N$ as the symplectic reflection algebra of $\N$ \cite{etingof2002symplectic}, it depends on a parameter $t$, $$A_{t}(W_{0},W):= \frac{T(V \oplus V^{*})\rtimes \N}{J}$$ where $J$ is an ideal of $\C \N$ generated by the relations $[x, x'] = 0$ for all $(x, x')\in V^{*} \times V^{*}$, $[y , y'] = 0$ for all $(y, y') \in V \times V$ and $[y, x] = x(y) + \sum_{H \in \mca_{0}}\frac{\alpha_{H}(y)x(v_{H})}{\alpha_{H}(v_{H})}\gamma_{H}$ for all $(x, y) \in V^{*} \times V$ and $v_{H} \in V$ such that $\C.v_{H} \oplus H = V$ and  $\C.v_{H}$ is $Fix_{W}(H)$ stable. The operator 
$$\gamma_{H}:= \sum_{j=0}^{m_{H}-1}m_{H}(k_{H,j+1}-k_{H,j})\epsilon_{H,j}$$ is a linear combination of primitive orthogonal idempotent of $\C \N$ where $m_{H}$ is the order of $Fix_{W}(H)$ and $(\epsilon_{H,j})_{j \in \{0, \ldots, m_{H}-1\}}$ is a system of primitives and orthogonal idempotents of $\C \N$ defined by $$\epsilon_{H,j}:= \frac{1}{|\N|}\sum_{w \in \N \setminus Id}det(w)^{j}w$$

In this article, we focus on the case of $t \neq 0$. Up to renormalization, we can reduce to $t=1$. Then the Cherednik algebra of the normalizer is denoted $A(W_{0},W)$.

According to the theorem 1.3 \cite{etingof2002symplectic}, we have an isomorphism of vector spaces $$A(W_{0},W) \simeq \C[V] \otimes \C \N \otimes \C [V^{*}]$$
We denote by $A(W_{0},W)_{reg}$ the localization by the multiplicative set defined by $\delta$, 

$$A(W_{0},W)_{reg}:= \C[X] \otimes_{\C[V]}A(W_{0},W)$$
\subsubsection{Dunkl operator} 
\begin{prp}(\cite{gobet2020hecke} prop. 2.6) The following differential 1-form is $W$-equivariant and integrable on $X$, 
$$\omega_{0} = \sum_{H \in \mca_{0}}a_{H}\frac{d\alpha_{H}}{\alpha_{H}} \in \Omega^{1}(X) \otimes \C W_{0}$$ where $a_{H} = \sum_{H \in \mca_{0}} m_{H}k_{H,j}\epsilon_{H,j}$.
\end{prp}
Let us define a family of Dunkl-Opdam operators by considering the covariant derivative of the connection $\nabla:= d + \omega_{0}$ on a trivial vector bundle over $X$.  
\begin{dfn}(Dunkl-Operator)  The covariant derivative of $\nabla$ along $y \in V$ is 
$$T_{y} := \partial_{y} + \sum_{H \in \mca_{0}} \frac{\alpha_{H}(y)}{\alpha_{H}}a_{H} \in \mcd(X) \rtimes \N$$ where $\mcd(X)$ stands for the algebra of differential operators of $X$ ( \cite{ginzburg1998lectures}).
\end{dfn}

\begin{prp}
\begin{enumerate}
    \item For all $(y, y') \in V \times V$, $[T_{y}, T_{y'}] = 0$.
    \item For all $g \in \N$ and for all $y \in V$, $gT_{y}g^{-1} = T_{g(y)}$.
\end{enumerate}
\end{prp}
\begin{proof}
For the commutativity, the proof is similar to that in \cite{etingof2002symplectic}.
\end{proof}
Then we define a faithful representation of $A(W_{0}, W)$.
\begin{thm}(The Dunkl embedding)
The morphism $$\fonctionb{\Phi}{A(W_{0},W)}{\mcd(X) \rtimes \N}{x \in V^{*}}{x}{y \in V}{T_{y}}{w \in \N}{w}$$ is an injective morphism of algebras.
\end{thm}
\begin{proof}
We define a filtration on $A(W_{0},W)$  by putting in degree $0$, $V^{*}$ and $\N$, and $V$ in degree $1$. 
We consider filtration by the order of differential operators on $\mcd(X) \rtimes \N$.

The morphism $\Phi$ is a morphism of graded algebras. So it induces a morphism on associated graded algebras $$gr(\Phi): gr(A(W_{0},W)) \rightarrow gr(\mcd(X)\rtimes \N)$$
We can prove that the composition $$\xymatrix{\C[V \oplus V^{*}] \rtimes \N \ar[r]^-{\simeq} & gr(A(W_{0},W)) \ar[r]^-{gr(\Phi)} & gr(\mcd(X) \rtimes \N) \ar[d]^{\simeq} \\ & & \C[X \oplus V^{*}] \rtimes \N }$$
is the identity on the homogenous components. Then  $gr(\Phi)$ is injective so is $\Phi$ .
\end{proof}

\begin{thm}
By localizing the previous morphism $\Phi$ becomes an isomorphism of algebras $$\Phi_{reg}: A(W_{0},W)_{reg} \simeq \mcd(X)\rtimes \N$$
\end{thm}
\begin{proof}
We just have to prove the surjectivity. The operator 

$\sum_{H \in \mca_{0}}\frac{\alpha_{H}(y)}{\alpha_{H}}a_{H}$ has no pole over $X$. Thus $y - \sum_{H \in \mca_{0}}\frac{\alpha_{H}(y)}{\alpha_{H}}a_{H}$ is an element of $A(W_{0}, W)_{reg}$. So the image of $\Phi_{reg}$ contains a system of generators of $\mcd(X) \rtimes \N$.
\end{proof}
\subsubsection{Category $\mco(W_{0}, W)$} We define  a category $\mco(W_{0}, W)$ similar to the category $\mco$ in \cite{ginzburg2003category}. We denoted by $eu_{0}:= \sum_{y \in \mathcal{B}}y^{*}.y - \sum_{H \in \mca}a_{H}$ the Euler element, where $\mathcal{B}$ is a basis of $V$.

The operator $\sum_{H \in \mca_{0}}a_{H}$ lies in the center of $\C \N$, so it acts by multiplication by a scalar $c_{E}$ on irreducible representations $E$ of $\N$. We define a partial order on simple $\C \N$-modules by $F < E$ if $c_{E} - c_{F} \in \mathbb{N}^{*}$.
\begin{prp}
\begin{enumerate}
    \item $[eu_{0},x] = x$ for all $x \in V^{*}$.
    \item $[eu_{0}, y] = y$ for all $y \in V$.
    \item $[eu,w] = 0$ for all $w \in \N$.
\end{enumerate}
\end{prp}
The operator $eu_{0}$ induces an inner graduation on $A(W_{0},W)$ defined by $A(W_{0},W)^{i} = \{ a \in A(W_{0},W) | [eu_{0}, a] = ia\}$ for all $i \in \mathbb{Z}$. The standard objects are $ \Delta(E):= Ind_{\C[V^{*}] \otimes \C \N}^{A(W_{0},W)}E$ where $E$ is a simple $\C \N$-module. The standard objects satisfy the same  properties as standard objects in the rational case, $\Delta(E)$ admits a simple head $L(E)$ and each simple head admits a projective cover $P(E)$ and each projective cover admit a filtration by $A(W_{0},W)$-submodules $P_{i}$ of $P(E)$ with successive quotient $P_{i}/P_{i-1}$ isomorphic to a standard object $\Delta(E_{i})$ such that $E_{i} < E$. We denote by $\Lambda$ the set of all simple objects. The category $\mco(W_{0}, W)$ is the full subcategory of finitely generated $A(W_{0}, W)$-module  satisfying the following  conditions: \begin{enumerate}
    \item The action of $\C[V^{*}]$ on $M \in \mco(W_{0}, W)$ is locally nilpotent.
    \item $M$ is isomorphic to the direct sum of its generalized $eu_{0}$-eigen space, $$M \simeq \bigoplus_{\alpha \in \mathbb{C}}\mathcal{W}_{\alpha}(M), \text{ with } \mathcal{W}_{\alpha}(M):= \{ m \in M | \exists N > 0, (eu_{0} - \alpha)^{N}.m = 0\}$$
\end{enumerate}
The standard objects $\Delta(E)$ are elements of $\mco(W_{0},W)$. The set $\Lambda$ of all simple objects is a set of complete simple objects of $\mco(W_{0},W)$. According to \cite{ginzburg2003category} theorem 2.19,  we can assert that the triple $(\mco(W_{0}, W), \Lambda, <)$ is a highest weight category in the sense of \cite{cline1988finite}. Therefore, it exists a quasi-hereditary cover, so there exists a finite dimensional $\C$-algebra $R$ such that $\mco(W_{0},W)$ is equivalent to the category of finitely generated modules over the algebra $R$.
 
 \subsection{The functor \texorpdfstring{$KZ_{0}$}.}

 We have a localization functor $\text{Loc}$ from the category of $A(W_{0}, W)$-modules to the category of $A(W_{0},W)_{reg}\text{-modules}$ which sends $M \in A(W_{0},W)\text{-mod}$ to $M_{reg}:= A(W_{0},W)_{reg} \otimes_{A(W_{0},W)}M$.
 
 Since $\C[X]$ is flat over $\C[V]$, this functor is exact. Let $M$ be a $A(W_{0},W)$-module, we  defined $M_{tor}$ as $\{m \in M | \exists N>0, \, \delta^{N}.m=0\}$. The full subcategory of $A(W_{0},W)$-modules such that $M = M_{tor}$ is called the category of   $A(W_{0},W)$-modules with torsion with respect of $\delta$ and it is denoted  $(A(W_{0},W)\text{-mod})_{tor}$. Thus, for all $M \in (A(W_{0},W)\text{-mod})_{tor}$, $\text{Loc}(M) = 0$. Then the localization functor factorized through the quotient category $\frac{A(W_{0},W)\text{-mod}} {(A(W_{0},W)\text{-mod})_{tor}}$. It induces a faithful functor $$\frac{A(W_{0},W)\text{-mod}} {(A(W_{0},W)\text{-mod})_{tor}} \rightarrow A(W_{0},W)_{reg}\text{-mod}$$
 We have a restriction functor $Res:A(W_{0},W)_{reg}\text{-mod}\rightarrow A(W_{0},W)\text{-mod}$ sending $M \rightarrow \C[X]\otimes_{\C[V]}M$ which is right adjoint to the localization functor. The functor $\text{Loc}$ induces an equivalence of categories $\frac{A(W_{0},W)\text{-mod}} {(A(W_{0},W)\text{-mod})_{tor}} \rightarrow A(W_{0},W)_{reg}\text{-mod}$ \cite{gabriel1962categories}.
 
 \medskip

 Let $\mco(W_{0},W)_{tor}:= \mco(W_{0},W) \cap (A(W_{0},W)\text{-mod})_{tor}$. For all objects with $\delta$-torsion i.e $M \in \mco(W_{0},W)_{tor}$, $\text{Loc}(M) = 0$. The functor  $\text{Loc}$ induces a fully faithful functor 
 $$\frac{\mco(W_{0},W)}{\mco(W_{0},W)_{tor}}\rightarrow A(W_{0},W)_{reg}\text{-mod}$$
 which is the composition of the functors 
 $$\frac{\mco(W_{0},W)}{\mco(W_{0},W)_{tor}}\rightarrow \frac{A(W_{0},W)\text{-mod}}{(A(W_{0},W)\text{-mod})_{tor}}$$
 and 
 $$\frac{A(W_{0},W)\text{-mod}}{(A(W_{0},W)\text{-mod})_{tor}} \rightarrow A(W_{0},W)_{reg}\text{-modules}$$

The first is fully faithful (lemma 3.3 \cite{rouquier2010derived}) and the second is the previous equivalence of categories. 

The Dunkl embedding induces an equivalence of categories between the category of $A(W_{0},W)\text{-modules}$ and  the category of $\mcd(X)\rtimes \N\text{-modules}$. 

Let $e:= \frac{1}{|\N|}\sum_{g \in \N}g$, this is a central idempotent of $\C \N$. The category of $\mcd(X)\rtimes \C \N \text{-modules}$ is equivalent to the category $e.(\mcd(X)\rtimes \C \N).e\text{-modules}$. Since we have $e.(\mcd(X)\rtimes \C \N).e \simeq \mcd(X)^{\N}$, the category of $\mcd(X)\rtimes \N$-modules is equivalent to the category of $\mcd(X)^{\N}\text{-modules}$.

\medskip

However, $\N$ is not necessarily a complex reflection group, so we do not have a priori an isomorphism between the algebra of differential operators invariants  by the action of $\N$ and the algebra of differential operators on $X/\N$. We use a general result from algebraic geometry proved by Halland and Cannigs in \cite{cannings1994differential} theorem 3.7. Let us recall this result. Let $k$ be an algebraically closed field. Let $A$ be a  reduced, finitely generated $k$-algebra. Let $X:= Spec(A)$. Let $G$ be a finite group acting on $A$ by automorphisms, then it acts on $X$. The algebra $A^{G}$ is  also reduced and finitely generated. We denote by $X/G:= Spec(A^{G})$. We denote by $\Phi: A^{G} \rightarrow A$ and $\Phi': Spec(A) \twoheadrightarrow Spec(A^{G})$, we get a morphism of sheaves $\Phi^{\sharp}: \mco_{Spec(A^{G})} \rightarrow \Phi_{*}'\mco_{Spec(A)}$. Let $p$ be a point in $X$ and $I^{A}_{G}(p):=\{ g \in G | (f-g.f) \in p, \, \fa f \in A\}$ the inertia group of $p$.

Let $\widetilde{V}:= \{ p \in X | I(p) = 1\}$. We said $G$ act generically  without inertia if $G$ acts trivially and if for all generic points of $X$, $I(p) = 1$. We can assume $G$ acts generically without  inertia. If $\widetilde{V} = X$, we stated that $G$ acts without fixed points.

\begin{thm}(\cite{cannings1994differential} theorem 3.7)
\begin{enumerate}
    \item If $G$ acts without fixed points then $\mcd(A^{G}) \simeq \mcd(A)^{G}$
\end{enumerate}
\end{thm}
Let us apply this theorem to the case of $G = \N$ and $A = \C[X]$. Since $\N$ acts on $X$ without fixed points,  we get the desired isomorphism $\mcd(X)^{\N} \simeq \mcd(X/\N)$. 

\medskip

Let us investigate the structure of $\mcd(X/\N)\text{-module}$. Let $E$ be a simple $\C \N$-modules. We have  $\Delta(E)_{reg} \simeq \C[X]\otimes E$, it is a free $\C[X]$-module of rank $dim(E)$. Thus, it corresponds to a trivial vector bundle over $X$ of fiber $E$. The structure of $\mcd(X) \rtimes \N$-module is given by the action of $x \in V^{*}$, $T_{y}$ with $y \in V$ and $w \in \N$ on $\Delta(E)_{reg}$. The element $w$ acts diagonally on $\Delta(E)_{reg}$, $x$ acts by multiplication on the left and $T_{y}(P \otimes v) = \partial_{y}(P \otimes v) + \sum_{H \in \mca_{0}}\frac{\alpha_{H}(y)}{\alpha_{H}}a_{H}(P \otimes v)$ where $P \in \C[X]$ and $v \in E$. Since $y.v = 0$ for all $y \in V$ and $v \in E$ then $T_{y}(v)=0$ so $\partial_{y}(v) = - \sum_{H \in \mca_{0}}\frac{\alpha_{H}(y)}{\alpha_{H}}a_{H}v$. Therefore, we get 
$$\partial_{y}(P \otimes v) = \partial_{y}P \otimes v + \sum_{H \in \mca_{0}}\frac{\alpha_{H}(y)}{\alpha_{H}} \sum_{j=0}^{m_{H}-1}m_{H}k_{H,j} P \otimes \epsilon_{H,j}.v$$
This formula defines a covariant derivative $\nabla^{0}_{y}:= \partial_{y}$. The associated connection is $\nabla^{0}:= d \otimes Id + \sum_{H \in \mca_{0}} \frac{d \alpha_{H}}{\alpha_{H}}(\sum_{j=0}^{m_{H}-1}m_{H}k_{H,j}(Id \otimes \epsilon_{H,j}))$.
\begin{prp}
The algebraic connection $\nabla^{0}$ is flat, $\N$-equivariant with regular singularities on $V$.
\end{prp}
\begin{proof}
According to the proposition 2.1 we obtain the flatness of the connection.

  Since for all $\alpha \in V^{*}$ and for all $g \in \N$, $g.d\alpha = d(g.\alpha)$ and  the operator $\sum_{H \in \mca_{0}}\frac{d \alpha_{H}}{\alpha_{H}}a_{H}'$ is $\N$-equivariant if and only if for all $H \in \mca_{0}$ and for all  $w \in \N$, $a_{w(H)}' = w a_{H}' w^{-1}$,  the connection $\nabla^{0}$ is $\N$-equivariant.

\medskip

 Let us  prove the claim about the singularities. Let us follow \cite{deligne2006equations},and try to apply the Deligne’s regularity criterion for an integrable connection
$\nabla$ on a smooth complex algebraic variety $X$. It states that $\nabla$ is regular along the irreducible divisor at infinity in some fixed normal compactification
of $X$ if and only if the restriction of $\nabla$ to every smooth curve on $X$ is
Fuchsian .

\medskip

Let $j:V \lra \mathbb{P}(V \oplus \C)=:Y$ be a compactificaction $\N$-equivariant of $V$. We consider  $M_{Y} := \mco_{Y} \otimes \Delta(E)$ a sheaf of $\mco_{Y}$-module.

\medskip
 
  Let $H \in \mca_{0}$ and  $x_{H} \in H$. We know that there is a vector $v_{H} \in V$ such that the line generated by  $v_{H}$ is  supplementary to $H$ and $W_{H}$-stable.
Let us consider the affine line directed by $L_{H}$ passing  through $x_{H}$ denoted $\Tilde{L_{H}}:= x_{H} + \C.v_{H}$.
 
Let $D_{H}^{*} =\{x_{H} + z.v_{H} | 0 < |z| < 2 \}$ such that $D_{H}^{*} \subset X$.
 Then $\alpha_{H}(x_{H} + z.v_{H}) = \alpha_{H}(x_{H}) + z.\alpha_{H}(v_{H}) = 0 + z.\alpha_{H}(v_{H})$.
 For $H \neq H'$, $\alpha_{H'}(x_{H} + z.v_{H}) \neq 0$ on $D_{H}^{*}$.
 We restricted the connection to the affine line, and we got
 \label{16}$$\nabla_{|_{\Tilde{L_{H}}}} = d \otimes Id - \frac{dz}{z}\sum_{j=0}^{m_{H}-1}m_{H}k_{H,j}Id \otimes \epsilon_{H,j} - \sum_{H' \neq H}\frac{\alpha_{H'}(v_{H})dz}{\alpha_{H'}(x_{H} + z.v_{H})} \sum_{j=0}^{m_{H}-1}m_{H} k_{H,j} Id \otimes \epsilon_{H,j}$$
It is deduced that on each $H$ the singularities are regulars.

\medskip

Now let us prove the regularity at infinity. We start by defining a change of chart.
Let $\phi: V \lra \C$ be a linear form of $V$. Let $\Check{\phi}$ such that $\phi(\Check{\phi}) = 1$. We define a chart of $Y$, noted $V_{\phi}$ by $V_{\phi} = \underbrace{\{ v \in V | \phi(v) = 1 \}}_{=:U_{\phi}} \oplus \C \subset V \oplus \C$. 
 Let  $(e_{1}, \ldots , e_{n})$ be the canonical basis of $V$ and $(e_{1}^{*}, \ldots , e_{n}^{*})$ its dual basis. The elements $V_{e_{i}^{*}}$ are the classical affine carts ($[x_{1}: \ldots : 1 : \ldots : x_{n+1}]$).
 We consider the pullback on $V_{\phi}$ of the $1$-form of connection. It is equivalent to a replacement of variables $(v,t) \lra (\frac{v}{t}, 1)$.
We get the formula:$$ d \otimes Id +\sum_{H \in \mca}\underbrace{a_{H}(\frac{-1}{t^{2}} \alpha_{H_{|_{U_{\phi}}}}(v))dt}_{=B} + \underbrace{\frac{\frac{1}{t} (d\alpha_{H})_{|_{U_{\phi}}}(v))}{(\frac{1}{t}\alpha_{H}(v))}}_{=A}  $$
Part  $A$ is regular by the previous step and does not depend on $t$. Part  $B$ is regular in $t = \infty$ after the change of variable $(t' = 1/t)$. Therefore, $\nabla^{0}$ is regular.
\end{proof}
Every standard object is endowed with a flat, $\N$-equivariant, connection with regular singularities on $V$. So every object of the category $\mco(W_{0},W)$  is endowed with a connection with regular singularities over $V$. We can apply the Riemann-Hilbert-Deligne equivalence, we get a $ \C \pi_{1}(X/\N)\text{-mod}_{f.d}$.
\medskip

If we compose this functor with the previous construction we get the following functor $\mco(W_{0},W) \rightarrow \C \pi_{1}(X/\N)\text{-mod}_{f.d}$, $M \rightarrow (((M_{reg})^{\N}))^{an})^{\nabla^{0}}$ where $(-)^{an}$ is the analytisation functor and $(-)^{\nabla^{0}}$ the horizontal section functor of $\nabla^{0}$.

\medskip

According to \cite{gobet2020hecke} proposition 2.6 the monodromy representation factorizes through $H(W_{0},W)$. Finally, we get a functor $$KZ_{0}: \mco(W_{0},W) \rightarrow H(W_{0},W)\text{mod}_{f.d}$$  such that $M \rightarrow (((M_{reg})^{\N})^{an})^{\nabla^{0}}$ which induces a faithful functor $\underline{KZ_{0}}: \frac{\mco(W_{0},W)}{\mco(W_{0},W)_{tor}}\rightarrow H(W_{0},W)\text{-mod}_{f.d}$. 
\medskip

The quasi-hereditary cover of $\mco(W_{0},W)$ and the exactness of $KZ_{0}$ allow us to use a kind of Watt's theorem. The claim of the following general proposition has been communicated to us by R.Rouquier.
\begin{prp}\label{thmrepresentabilite}
Let $A$ and $B$ two $k$-algebras of finite dimensions. Let $F$ be an exact functor  from the category of finitely generated $A$-module to the category of finitely generated $B$-modules. Then $F$ is isomorphic to the functor $Hom_{A}(Hom_{A^{op}}(F(A),A),-)$.
\end{prp}
\begin{proof}
What follows is a mere outline of the proof: firstly, based on the proposition 4.4.b \cite{auslander1997representation} the functor $Hom_{A}(F(A), A) \otimes_{A}-$ is isomorphic to the functor $Hom_{A}(F(A),-)$. Secondly, based on the corollary 5.47 \cite{rotman2008introduction} the functor $F(A)\otimes_{A}-$ is isomorphic to the functor $F(-)$. Thirdly, based on proposition 4.3.b \cite{ariki1995representation}, $Hom_{A}(Hom_{A^{op}}(F(A),A),A)$ is isomorphic to $F(A)$.
\end{proof}
Let us apply this theorem to our example, the functor $KZ_{0}$ is representable by a projective object of $\mco(W_{0},W)$ called $P_{KZ_{0}}$. 
\begin{thm}
The morphism of algebra $\Phi: H(W_{0},W) \rightarrow End_{\mco}(P_{KZ_{0}})^{op}$ is actually an isomorphism.
\end{thm}
\begin{proof}
Let us follow the arguments proposed in \cite{bellamy2012symplectic} part 4.6.
 \begin{lem}(\cite{bellamy2012symplectic} lemma 4.6.4)\label{lemtech}
Let $\mca$ be an abelian category. We assume $\mca$ is Artinian. Let $\mcb$ be an abelian full subcategory of $\mca$, closed under sub-object and quotient.
 Then the functor $\mcf: \mcb \lra \mca$ has a left adjoint  $\mcg: \mca \lra \mcb$ sending an object $M \in Obj(\mca)$ on its largest quotient in $\mcb$ and the co-unit $\eta: Id_{\mca} \lra \mcf \circ \mcg$ induces a family of surjective morphisms for all objects in $\mca$.
 \end{lem}
 Let us apply this lemma to $\mca = H(W_{0},W)\text{-mod}_{f.d}$ and $\mcb = Im(KZ_{0})$. This proves the surjectivity of $\Phi$. 
 
 Since $$P_{KZ_{0}} = \bigoplus_{E \in Irr(\N)}dim(KZ_{0}(L(E)))P(E)$$ we can calculate explicitly $dim(End_{\mco(W_{0},W)}(P_{KZ_{0}})^{op})$.
 
 We have $dim(End_{\mco(W_{0},W)}(P_{KZ_{0}})^{op}) = dim(H(W_{0},W))$. This implies that $\Phi$ is an isomorphism.
\end{proof}
\begin{thm}
The functor $\underline{KZ_{0}}$ is an equivalence of categories.
\end{thm}
\begin{proof}
The only thing which remains to prove is the essential surjectivity. The category $Im(KZ_{0})$ is a full subcategory of $H(W_{0},W)\text{-mod}_{f.d}$, closed under quotients, sub-object and direct sum. Since $H(W_{0},W)$ is isomorphic to $End_{\mco(W_{0},W)}(P_{KZ_{0}})^{op}$,  $H(W_{0},W)$ is an object of $Im(KZ_{0})$, so $H(W_{0},W)^{m}$ is an object of $Im(KZ_{0})$. Then every $H(W_{0},W)\text{-modules}$ finitely generated lies into $Im(KZ_{0})$. Therefore, $\underline{KZ_{0}}$ is essentially surjective.
\end{proof}
\subsubsection{Forgetting the ambiant group $W$}
We can provide a related result involving only $W_{0}$, and not the ambient group $W$. This is done in our article \cite{fallet2022cherednik}. The general setting is as follows.
Let $G$ be a finite subgroup of $GL(V)$. Let $G_{0}$ be a normal subgroup of $G$ generated by reflexions. Let $\mcr_{0}$ be the set of reflexions of $G_{0}$ and $\mca_{0}$ the arrangement of reflecting hyperplanes of $G_{0}$. The first goal is to build up a Hecke algebra for $G$ from the Hecke algebra of $G_{0}$ generalizing $H(W_{0},W)$ for $G = \N$.

Let $X^{+}$ be the subspace of $V$ on which $G$ acts freely and let $X_{0}$ be the subspace of $V$ on which $G_{0}$ acts freely. The manifold $X_{0} \setminus X^{+}$ is of codimension $>2$ then $\pi_{1}(X^{+}) \simeq \pi_{1}(X_{0})$ \cite{godbillon1971elements} theorem 2.3. 
We get two short exact sequences.

$$\xymatrix{1 \ar[r] & \pi_{1}(X^{+}) \ar[r] & \pi_{1}(X^{+}/G_{0}) \ar[r] & G_{0} \ar[r] & 1 \\
1 \ar[r] \ar[u]^{=} & \pi_{1}(X_{0}) \ar[r] \ar[u]^{\simeq} & \pi_{1}(X_{0}/G_{0}) \ar[r] \ar[u] & G_{0} \ar[r] \ar[u]^{=} & 1 \ar[u]^{=}}$$

The exactness and the commutativity of the diagram  together imply $$\pi_{1}(X^{+}/G_{0}) \simeq \pi_{1}(X_{0}/G_{0})$$

The braid group $B_{0}$ of $G_{0}$ is a normal subgroup of $ B := \pi_{1}(X^{+}/G)$, we get a short exact sequence
$$\xymatrix{1 \ar[r] & B_{0}:= \pi_{1}(X_{0}/G_{0}) \ar[r] & \pi_{1}(X^{+}/G) \ar[r] & G/G_{0} \ar[r] & 1}$$

Let $I$ be the ideal of $\C B_{0}$ generated by the relations $\s_{H}^{m_{H}} = \sum_{k = 0}^{m_{H}-1}a_{H,k}\s_{H}^{k}$ for $\s_{H}$ a braided reflection associated to $H \in \mca_{0}$. Then the Hecke algebra of $G_{0}$ is the quotient $H_{0} := \frac{\C B_{0}}{I}$. According to the  now proven BMR freeness conjecture (see the references of \cite{gobet2020hecke} or its weaker version in Characteristic $0$ \cite{etingof2017proof}) it is an algebra finitely generated of dimension $|G_{0}|$. Let $I^{+} = \C B \otimes_{\C B_{0}} I$ be the ideal which define the Hecke algebra of $G$, $H(G):= \frac{\C B}{I^{+}} \simeq \C B \otimes_{\C B_{0}} H_{0}$ is of dimension $|G|$.

Let us make a link between this new algebra and the algebra $H(W_{0},W)$.
We defined $H(W_{0},W)$ as a quotient of the algebra $\C \Tilde{B}_{0}$. We defined $\Tilde{B}_{0}$ as the quotient of $\pi_{1}(X/ \N)$ by $K:= Ker(\pi_{1}(X) \rightarrow \pi_{1}(X_{0}))$. Since $X_{0} \setminus X^{+}$ has codimension $>2$ $$K =  Ker(\pi_{1}(X) \rightarrow \pi_{1}(X_{0})) \simeq Ker(\pi_{1}(X/\N) \rightarrow \pi_{1}(X^{+}/\N))$$  

And $\Tilde{B}_{0} \simeq \pi_{1}(X^{+}/ \N)$ is our group $\pi_{1}(X^{+}/G) =: B $. As a result, the algebra $H(W_{0},W)$ is the same as $H(G)$.

\vspace{0.5cm}

 Let us consider the category $\mco_{tor}^{0}$ the full subcategory of $\mco$ of module annihilated by a power of $\delta_{0}:= \prod_{H \in \mca_{0}}\alpha_{H}$. We have 
 \begin{thm}(\cite{fallet2022cherednik} theorem 5 )
$KZ_{0}$ is fully faithful and essentially surjective from the category $\frac{\mco}{\mco_{tor}^{0}}$ to the category  of finite dimension $H(G) \text{-modules}$.
\end{thm}
 A priori $\mco_{tor}$ and $\mco_{tor}^{0}$ are different. Actually, we can prove that these two categories are the same \cite{fallet2022cherednik}.
\section{The  rational Cherednik algebra of the pair \texorpdfstring{$(\mcl,W)$}- and the functor \texorpdfstring{$\widetilde{KZ}$}-.}
\subsection{The  rational Cherednik algebra of the pair \texorpdfstring{$(\mcl,W)$}.}
Let us denote the Cherednik algebra of the pair $(\mcl,W)$ by $A(\mcl,W)$. As a vector space $A(\mcl,W)$ is $\C[V]\otimes \C \mcl \otimes \C W \otimes \C [V^{*}]$. Let us define a product of algebra by adding relations between generators of the sub-algebras $\C[V]$, $\C[V^{*}]$, $\C \mcl \rtimes W$: $[x, x'] = 0$ for all $(x, x')\in V^{*} \times V^{*}$, $[y, y'] = 0$ for all $(y, y') \in V \times V$, $[e_{H}, e_{H'}] = 0$ for all $(H, H') \in \mcl \times \mcl$ and $$[y, x] = t. x(y) + \sum_{H \in \mca} \frac{\alpha_{H}(y)x(v_{H})}{\alpha_{H}(v_{H})}\gamma_{H}e_{H}$$ where $t \in \mathbb{C}$. Likewise to the case of $A(W_{0},W)$, we just consider the case $t \neq 0$. Up to renormalization, we can consider $t = 1$.
\subsubsection{Dunk operator}
In \cite{marin2018lattice} is introduced a differential 1-form, $W$-equivariant and integrable
$$\Tilde{\omega}:= \sum_{H \in \mca}\frac{d \alpha_{H}}{\alpha_{H}}a_{H}e_{H} \in \Omega^{1}(X) \otimes \C \mcl \otimes \C W$$
We define the Dunkl operator as the covariant derivative associated to the connection $\nabla := d  + \Tilde{\omega}$, in the direction of $y \in V$.
\begin{dfn}(Dunkl operator of the pair $(\mcl, W)$)
For all $y \in V$, $$\widetilde{T_{y}}:= \partial_{y} + \sum_{H \in \mca} \frac{\alpha_{H}(y)}{\alpha_{H}}a_{H}e_{H} \in (\mcd(X) \otimes \C \mcl)\rtimes W$$
\end{dfn}
This family of differential operators satisfied the following properties:
\begin{prp}
\begin{enumerate}
    \item For all $y, y' \in V \times V$, $[\widetilde{T_{y}}, \widetilde{T_{y'}}] = 0$.
    \item For all $y \in V$ and $g \in W$, $g.\widetilde{T_{y}}.g^{-1} = \widetilde{T_{g(y)}}$
\end{enumerate}
\end{prp}
\begin{proof}
\begin{enumerate}
    \item We can follow the method of \cite{etingof2010lecture} owing to the fact that the action of $(\mcd(X) \otimes \C^{c})\rtimes W$  on $\C(V)$ is faithful, with $c \in \mcl/W$. Let us prove a more general result. An outer automorphism of rings $R$ is an automorphism of $R$ which  is not inner. If each  non identity element of a group $G$ induces an outer automorphism of $R$, then $G$ is called a group of outer automorphisms. 
    \begin{prp}
    Let $R$ be a simple ring with unity. Let $G$ be a group of outer automorphisms. Let $X$ be a finite set. If $G$ acts transitively on $X$, then $(R \otimes \C^{X})\rtimes G$ is simple.
    \end{prp}
    \begin{proof}
    Let us start by proving the reducibility of the support of an element of $(R \otimes \C^{X})\rtimes G$.
    Let $J$ be a two-sided ideal of $(R \otimes \C^{X})\rtimes G$ not reduced to $\{0\}$. There exists $x \in J$, $x \neq 0$, of minimal support, $x = \sum_{g \in G}r_{g}.g = \sum_{g \in G, \lambda \in X} r_{\lambda,g}\epsilon_{\lambda}g$.
    Since $x$ is different from $0$, $(r_{g})_{g}$ are not all $0$. There exists $\lambda_{0} \in X$ such that $r_{\lambda_{0},g_{0}}.\epsilon_{\lambda_{0}} \neq 0$. So 
    \begin{align*}
        g_{0}^{-1}x &= \sum_{g \in G} g_{0}^{-1}r_{g}g = g_{0}^{-1}r_{g_{0}}g_{0} + \sum_{g \in G,  g \neq g_{0}}g_{0}^{-1}r_{g}g \\
                   &= r_{1}.1 + \sum_{g \in G,  g \neq g_{0}}g_{0}^{-1}r_{g}g
    \end{align*}
    Up to multiply $x$ by $g_{0}^{-1}$, we can assume  the existence of $\lambda_{0} \in X$ such that  $r_{1, \lambda_{0}}\epsilon_{\lambda_{0}} \neq 0$. The ideal $R r_{1, \lambda_{0}}R$ is a two-sided ideal of $R$ not reduced to $\{0\}$. Since $R$ is simple $R r_{1, \lambda_{0}}R = R$, so there exists $(x_{j}, y_{j}) \in R \times R$, such that $\sum_{j}x_{j}r_{1, \lambda_{0}}y_{j} = 1$. All depends on $\lambda_{0}$. Let us cancel this dependency. 
    
    We get $r_{1} = \sum_{\lambda \in X}r_{\lambda,1}\epsilon_{\lambda} = r_{\lambda_{0},1} \epsilon_{\lambda_{0}} + \sum_{\lambda \in X, \lambda \neq \lambda_{0}}r_{\lambda,1}\epsilon_{\lambda}$. 
    Let us consider $x'= x.\epsilon_{\lambda_{0}} \in J \setminus \{0\}$.
    Then \begin{align*}
        x' &= r_{\lambda_{0},1}\epsilon_{\lambda_{0}}\epsilon_{\lambda_{0}}.1 + \sum_{g \neq 1}r_{g}'.g \text{ where } r_{g}' = \sum_{\lambda \in X}r_{\lambda,g}\epsilon_{\lambda}\epsilon_{\lambda_{0}} = r_{\lambda_{0},g}\epsilon_{\lambda_{0}} \\
        &= r_{\lambda_{0},1}\epsilon_{\lambda_{0}} + \sum_{g \neq 1}r_{\lambda_{0},g}\epsilon_{\lambda_{0}}g
    \end{align*}
    but $\sum_{j}x_{j}x'y_{j} = \sum_{g \in G, j}x_{j}r_{g}'gy_{j} = \sum_{g \in G, j}x_{j}r_{g}' g(y_{j})g$. For $g = 1$, we have
    \begin{equation*}
         P_{1}':= \sum_{j}x_{j}r_{1}y_{j}= \sum_{j}x_{j}r_{\lambda_{0},1}\epsilon_{\lambda_{0}} y_{j}=\epsilon_{\lambda_{0}}
    \end{equation*}
    We can reduce $P_{1}'$ to $1$. Let 
    \begin{equation*}
        x" := \sum_{h \in G} hx'h = \sum_{h,g \in G}h(r_{g}') hgh^{-1} = \sum_{h \in G}h(r_{1}') + \sum_{h \in G, g \neq 1}hgh^{-1}
    \end{equation*}
    but $\sum_{h \in G}h(r_{1}') = \sum_{h \in G}h.\epsilon_{\lambda_{0}} = \sum_{h \in G}\epsilon_{h. \lambda_{0}} = cst . 1$.
    Finally, the support of $\sum_{j}x_{j}x"y_{j}$ is equal to the support of $x$ and 
    \begin{align*}
        [\sum_{j}x_{j}x"y_{j},x] &= \sum_{g \in G} P_{g}"gr - rP_{g}"g \\
                                 &= P_{1}"r - rP_{1}" + \sum_{g \in G, g \neq 1}(P_{g}"gr - rP_{g}"g)
    \end{align*}
    At least, one term misses in the support of $[x,r]$, so by the minimality of the support of $x$, $[x,r] = 0$. So we can reduce the support of an element. 
    
    Let $\Tilde{R}:= R \otimes \C^{X} = \bigoplus_{\lambda \in X}R.\epsilon_{\lambda}$. We have $[x,r] = \sum_{g \in G}(r_{g}g(r) - r.r_{g}).g =0$, due to the fact that $(R \otimes \C^{X})\rtimes G$ is a $\C G$-module, where the elements of $G$ form a basis, $r_{g}.g(r) = r.r_{g}$ for all $g \in G$ and $r \in R \otimes \C^{X}$. So $\Tilde{R}r_{g} = r_{g}\Tilde{R}$. 
    
    Let $g \in G$ such that $r_{g} \neq 0$, we can assume the existence  of $\mu_{0} \in X$ such that $r_{\mu_{0},g} \neq 0$. Then $r_{g}\epsilon_{g(\lambda)} = \sum_{\lambda \in X}r_{\mu,g}\epsilon_{\mu} \epsilon_{g(\lambda)} = r_{g(\lambda),g}\epsilon_{g(\lambda)}$ and 
    \begin{equation*}
        \epsilon_{\lambda}r_{g} = \sum_{\lambda \in X}\epsilon_{\lambda}r_{\mu,g}\epsilon_{\mu} = \sum_{\mu \in X}r_{\mu,g}\epsilon_{\lambda}\epsilon_{\mu} = r_{\lambda,g}\epsilon_{\lambda}
    \end{equation*}

    For $\lambda = \mu_{0}$, we get  $r_{g}\epsilon_{g(\mu_{0}} = r_{g(\mu_{0}),g}\epsilon_{\mu_{0}}=\epsilon_{\mu_{0}}r_{g}$; so $g(\mu_{0}) = \mu_{0}$. As a result, for all $\mu_{0}$ such that $r_{\mu_{0},g} \neq 0$ then  $g(\mu_{0}) = \mu_{0}$. Let $r$ be an element of $R\epsilon_{\mu_{0}}$, we can write $r = \rho \epsilon_{\mu_{0}}$ and $r_{\mu_{0},g} \neq 0$,
    \begin{align*}
        r_{g}g(r) &= r_{g}(g.\rho)\epsilon_{g.\mu_{0}} \\
                  &= \sum_{\mu \in X}r_{\mu,g} \epsilon_{\mu}g \rho \epsilon_{g(\mu_{0})} = \sum_{\mu \in X} r_{\mu,g}g.\rho \epsilon_{\mu}\epsilon_{g(\mu_{0})} \\
                  &= r_{g(\mu_{0}),g}g.\rho.\epsilon_{g(\mu_{0})} = g.\rho r_{\mu_{0},g}\epsilon_{\mu_{0}}
    \end{align*}
    and 
    \begin{equation*}
        \rho^{2}\mu_{0}r_{g} = \rho \epsilon_{\mu_{0}} \sum_{\mu \in X} r_{\mu,g}\epsilon_{\mu} = \sum_{\mu \in X} \rho r_{\mu,g} \epsilon_{\mu_{0}}\epsilon_{\mu} = \rho r_{\mu_{0},g}\epsilon_{\mu_{0}} = g.\rho r_{\mu_{0},g}\epsilon_{\mu_{0}}
    \end{equation*}
    Because of $r_{\mu_{0},g}g.\rho = \rho.r_{\mu_{0},g}$,  $r_{\mu_{0},g}R = R r_{\mu_{0},g}$ then $r_{\mu_{0},g}$ is two-sided invertible. The action of $g$ on $\rho$ would be by inner automorphism, which is absurd.
    \end{proof}
    Let us apply this result to $X = c \in \mcl/W$, $G = W$ and $R = \mcd(X)$.
    The algebra $(\mcd(X) \otimes \C^{c})\rtimes W$ is  simple, therefore $\bigoplus_{c \in \mcl/W} ( \mcd(X) \otimes \C^{c})\rtimes W$ is simple. But this algebra is $\mcd(X) \otimes \C \mcl) \rtimes W$. So the action of $(\mcd(X) \otimes \C \mcl)\rtimes W$ on $\C(V)$ is fully faithful. The rest of the proof is similar to the proof for the theorem 2.16 in \cite{etingof2010lecture}.
    \end{enumerate}
\end{proof}
\subsubsection{Dunkl embedding.}
We consider the filtration $(\mcf_{i})_{i \in \mathbb{Z}}$ on $(\mcd(X) \otimes \C \mcl) \rtimes W$ by the order of the differential operator.
Let us define a graduation on $(\C[X \oplus V^{*}] \otimes \C \mcl) \rtimes W$, by $$B_{i} = \{\sum\limits_{\substack{\alpha \in \mathbb{N}^{n}, \beta \in \mathbb{N}^{m}, \\ \lambda \in \mcl, g \in W, \\ |\beta| = i}} a_{\alpha, \beta, \lambda, g} x^{\alpha}y^{\beta}e_{\lambda}g\}$$ for all $i \in \mathbb{Z}$.
We get an isomorphism of algebras $\psi: grad_{\mcf}(\mcd(X) \otimes \C \mcl) \rtimes W) \simeq (\C[X \oplus V^{*}] \otimes \C \mcl) \rtimes W$. 
Let us define a filtration $(F_{i})_{i \in \mathbb{Z}}$ on $A(\mcl,W)$ by putting in degree $0$, $\mcl$, $V^{*}$, $W$ and $V$ in degree $1$.
Let us define a graduation on $(\C[V \oplus V^{*}] \otimes \C \mcl)\rtimes W$ by $$A_{i} := \{ \sum\limits_{\substack{\alpha \in \mathbb{N}^{n}, \beta \in \mathbb{N}^{m}, \\ g \in W, \lambda \in \mcl, \\ | \beta| = i}}a_{\alpha, \beta, \lambda, g} x^{\alpha}y^{\beta}e_{\lambda}g \}$$
We get an isomorphism of algebras $\Phi:grad_{F}(A(\mcl,W)) \simeq \C[V \oplus V^{*}] \otimes \C \mcl)\rtimes W$.

\begin{thm}
The application $$\fonctionb{\Phi}{A(\mcl,W)}{(\mcd(X) \otimes \C \mcl) \rtimes W}{x \in V^{*}}{x}{y \in V}{\Tilde{T_{y}}}{g \in W}{g}$$ is an injective morphism of algebras. After localization, $\Phi$ becomes an isomorphism of algebras from $A(\mcl,W)_{reg}$ to $(\mcd(X) \otimes \C \mcl) \rtimes W$, denoted $\Phi_{reg}$.
\end{thm}
\begin{proof}
We observe that $\Phi$ is a morphism of filtered algebras because of $\Phi(F_{i}) \subset \mcf_{i}$. Therefore, $\Phi$ induces a morphism of algebras between graded algebras denoted by $gr(\Phi)$. Then we can  prove that $gr(\Phi)$ is injective by considering the following  composition $$\xymatrix{(\C[V \oplus V^{*}] \otimes \C \mcl)\rtimes W \ar[r]& \text{grad}_{F}(A(\mcl,W)) \ar[r]& \text{grad}_{\mcf}(\mcd(X) \otimes \C \mcl)\rtimes W) \ar[d] \\ & & (\C[X \oplus V^{*}] \otimes \C \mcl)\rtimes W)}$$
This is the identity. Therefore, $gr(\Phi)$ is injective, which implies that the morphism $\Phi$ is injective.

Since $\sum_{H \in \mca}\frac{\alpha_{H}(y)}{\alpha_{H}}a_{H}e_{H}$ is well defined on $A(\mcl, W)_{reg}$, the image of $\Phi_{reg}$ contains a system of generators of $(\mcd(X) \otimes \C \mcl) \rtimes W$.
\end{proof}

\subsubsection{Category $\mco(\mcl,W)$.}
Let $\Tilde{eu}:= \sum_{y \in \mathcal{B}}y^{*}.y - \sum_{H \in \mca}a_{H}e_{H}$ be the Euler element, where $\mcb$ is a basis of $V$.

\begin{prp}
$[\Tilde{eu},x] = x$ for all $x \in V^{*}$, $[\Tilde{eu},y] = -y$ for all $y \in V$, $[\Tilde{eu}, e_{H}] = 0$ for all $e_{H} \in \C \mcl$, $[\Tilde{eu}, w] = 0$ for all $w \in W$.
\end{prp}
Therefore, $\Tilde{eu}$ induces an inner graduation on $A(\mcl,W)$, $$A(\mcl,W)^{i}:=\{ a \in A(\mcl,W) | [\Tilde{eu},a] =ia\}$$ for all $i \in  \mathbb{Z}$.

The element $\sum_{H \in \mca}a_{H}e_{H}$ belongs to $Z(\C \mcl \rtimes W)$. Thus, it acts by scalar multiplication   on simple $(\C \mcl \rtimes W)$-modules. We denote by $c_{E}$ the associated scalar where $E$ is a simple $(\C \mcl \rtimes W)$-module. 

\medskip

Let us define a partial order on simple $\C \mcl \rtimes W$-modules by: $E' < E$ if $c_{E} - c_{E'} \in \mathbb{N}^{*}$. Let us define the standard objects associated to a simple $\C \mcl \rtimes W$-module, by $\Delta(E):= Ind_{\C[V^{*}] \otimes (\C \mcl \rtimes W)}^{A(\mcl,W)}E$. This object admits a simple head $L(E)$ and each simple head admits a projective cover $P(E)$. Each projective cover $P(E)$ admits a standard filtration. We denote by $\Lambda$ the set of all simple heads. 

\medskip

The category $\mco(\mcl,W)$ is the full subcategory of $A(\mcl,W)$-module finitely generated, locally nilpotent  for the action of $\C[V^{*}]$ and isomorphic to the direct sum of their generalized $\Tilde{eu}$-eigen spaces. The triple $(\mco(\mcl,W), \Lambda, <)$ is a highest weight category in the sense of \cite{cline1988finite}. Then $\mco(\mcl,W)$ admits a quasi-hereditary cover, this means that there exists a finite dimensional $\C$-algebra $\Tilde{R}$ such that  $\mco(\mcl,W)$ is equivalent to the category of finitely generated $\Tilde{R}$-modules.
\subsection{The functor \texorpdfstring{$\widetilde{KZ}$}.}
Let $M$ be a $A(\mcl,W)$-module. Let $M_{tor}:=\{ m \in M | \exists n >0 \, \delta^{n}.m = 0\}$ and $(A(\mcl,W)\text{-module})_{tor} = \{ M \in A(\mcl,W)\text{-mod} | M_{tor} = M\}$. This is a Serre subcategory of $A(\mcl,W)$-module. 

Let $\text{Loc}:A(\mcl,W)\text{-modules} \rightarrow A(\mcl,W)_{reg}\text{-modules}$ the localization functor, $\text{Loc}(M) = M_{reg}:= \C[X] \otimes_{\C[V]}M$, this is an exact functor, because $\C[X]$ is a flat $\C[V]$-module and $\text{Loc}(M) = 0$ if and only if $M \in (A(\mcl,W)\text{-mod})_{tor}$. Therefore, the functor $\text{Loc}$ factorizes through the quotient category $\frac{A(\mcl,W)\text{-mod}}{(A(\mcl,W)\text{-mod})_{tor}}$ and then the induced quotient functor $\frac{A(\mcl,W)\text{-mod}}{(A(\mcl,W)\text{-mod})_{tor}} \rightarrow A(\mcl,W)_{reg}\text{-mod}$ is an equivalence of categories. 

We can restrict this functor to the category $\mco(\mcl,W)$. Let us introduce the category $\mco(\mcl,W)_{tor}:= \mco(\mcl,W) \cap (A(\mcl,W)\text{-mod})_{tor}$. This is a Serre subcategory of $\mco(\mcl,W)$. The lemma 3.3 \cite{rouquier2010derived} implies that $\frac{\mco(\mcl,W)}{\mco(\mcl,W)_{tor}}\rightarrow \frac{A(\mcl,W)\text{-mod}}{(A(\mcl,W)\text{-mod})_{tor}}$ is a fully faithful functor.
By composing this with the previous equivalence of categories, we get a fully faithful functor $\frac{\mco(\mcl,W)}{\mco(\mcl,W)_{tor}}\rightarrow A(\mcl,W)_{reg}\text{-mod}$. Then  we apply the Dunkl embedding. 

Let us figure out the structure of $(\mcd(X) \otimes \C \mcl) \rtimes W$-module, on a standard object $\Delta(E)$, $E \in Irr(\C \mcl \rtimes W)$. The localization $\Delta(E)_{reg}$ of $\Delta(E)$ correspond to a trivial vector bundle over $X$ of dimension $dim(E)$. We can endow this vector bundle  with a connection by considering the action of $\Tilde{T_{y}}$ on an element $P \otimes v \in \Delta(E)_{reg}$. We get the formula $$\nabla_{y}(P\otimes v):= \partial_{y}P \otimes v + \sum_{H \in \mca} \frac{\alpha_{H}(y)}{\alpha_{H}} \sum_{j=0}^{m_{H}-1}m_{H}k_{H,j} P\otimes \epsilon_{H,j}v$$
\begin{prp}
$\nabla_{y}$ is an algebraic, flat and $W$-equivariant connection with regular singularities on $V$.
\end{prp}
\begin{proof}
The proof is similar to that proposed for the $A(W_{0},W)$ case.
\end{proof}
Since the connection $\nabla_{y}$ has regular singularities, we can apply the Riemann-Hilbet-Deligne correspondence and we get  a $\C \pi_{1}(X/W) \ltimes \C \mcl\text{-mod}_{f.d}$. The category of connection with regular singularities is thick and for all simple $\C \mcl \rtimes W$-modules $E$ we get a short exact sequence $$\xymatrix{\Delta(E) \ar[r] & L(E) \ar[r] & 0}$$ So $L(E)$ is endowed with a connection with regular singularities. Every object of $\mco(\mcl,W)$ admits a finite Jordan-Hölder series, so all objects of $\mco(\mcl,W)$ can be endowed with a connection with regular singularities.

According to  proposition 5.6 and 5.7 \cite{gobet2020hecke} this monodromy action factorizes through $\mcc(\mcl,W)$. We obtain an exact functor $$\fonction{\widetilde{KZ}}{\mco(\mcl,W)}{\mcc(\mcl,W)\text{-mod}_{f.d}}{M}{(((M_{reg})^{W})^{an})^{\nabla}}$$
According to \ref{thmrepresentabilite} the functor $\widetilde{KZ}$ is representable by a projective object in $\mco(\mcl,W)$, denoted by $P_{\widetilde{KZ}}$. The image of the functor $\widetilde{KZ}$ is a full abelian subcategory of the category of $\mcc(\mcl,W)$-modules finitely generated, closed under quotient, sub-objects and direct sum. We get the result:
\begin{prp}
The morphism $\Phi:\mcc(\mcl,W) \rightarrow End_{\mco(\mcl,W)}(P_{\widetilde{KZ}})$ is an isomorphism of algebras.
\end{prp}
\begin{thm}
The functor $\widetilde{KZ}$ is essentially surjective. Hence, the induced functor $\underline{\widetilde{KZ}}:\frac{\mco(\mcl,W)}{\mco(\mcl,W)_{tor}} \rightarrow \mcc(\mcl,W)\text{-modules}_{f.d}$ is an equivalence of categories.
\end{thm}
The proofs of these results are similar to those in the previous section.

\section*{Acknowledgements:}
These results are part of my PhD-Thesis at University Picardie Jules Verne under the supervision of Prof. Ivan Marin. I would like to thank Cedric Bonnafé, which suggested forgetting the ambient group W.


\bibliographystyle{alpha} 
  \bibliography{ms.bib}
\end{document}